\documentclass[11pt]{article}
\usepackage[utf8]{inputenc}
\usepackage{amsmath,amssymb,epsfig,bbm}
\usepackage{stmaryrd,mathabx}
\usepackage{comment}
\usepackage{color}
\usepackage[T1]{fontenc}

\usepackage[textsize=small]{todonotes}
\usepackage{enumitem}
\usepackage{varwidth}
\setlist{nolistsep}
\usepackage{hyperref}

\newtheorem{defi}{Definition}[section]
\newtheorem{prop}[defi]{Proposition}
\newtheorem{theo}[defi]{Theorem}
\newtheorem{conj}[defi]{Conjecture}
\newtheorem{lemm}[defi]{Lemma}
\newtheorem{coro}[defi]{Corollary}
\newtheorem{rema}[defi]{Remark}
\newtheorem{exem}[defi]{Example}
\newtheorem{exems}[defi]{Examples}

\newcommand{\bdefi}{\begin{defi}}
\newcommand{\edefi}{\end{defi}}
\newcommand{\bprop}{\begin{prop}}
\newcommand{\eprop}{\end{prop}}
\newcommand{\btheo}{\begin{theo}}
\newcommand{\etheo}{\end{theo}}
\newcommand{\blemm}{\begin{lemm}}
\newcommand{\brema}{\begin{rema}}
\newcommand{\erema}{\end{rema}}
\newcommand{\bexer}{\begin{exem}}
\newcommand{\eexer}{\end{exem}}
\newcommand{\bexems}{\begin{exems}}
\newcommand{\eexems}{\end{exems}}
\newcommand{\bconj}{\begin{conj}}
\newcommand{\econj}{\end{conj}}
\newcommand{\elemm}{\end{lemm}}
\newcommand{\bcoro}{\begin{coro}}
\newcommand{\ecoro}{\end{coro}}
\newcommand{\dem}{\noindent{\bf Proof. }}

\usepackage{mathrsfs}
\renewcommand\mathcal{\mathscr}

\renewcommand{\H}{{\cal H}}

\newcommand{\M}{{\cal M}}
\newcommand{\N}{{\cal N}}
\newcommand{\OOO}{{\cal O}}

\newcommand{\maths}[1]{{\mathbb #1}}  

\newcommand{\HH}{\maths{H}}
\newcommand{\II}{\maths{I}}
\newcommand{\JJ}{\maths{J}}
\newcommand{\KK}{\maths{K}}

\newcommand{\NN}{\maths{N}}

\newcommand{\PP}{\maths{P}}
\newcommand{\QQ}{\maths{Q}}
\newcommand{\RR}{\maths{R}}
\newcommand{\SSS}{\maths{S}}

\newcommand{\weakstar}{\overset{*}\rightharpoonup}
\newcommand{\ra}{\rightarrow}
\newcommand{\bs}{\backslash}

\newcommand{\ga}{\gamma}
\newcommand{\Ga}{\Gamma}

\newcommand{\cqfd}{\hfill$\Box$}

\newcommand{\bigO}{\operatorname{O}}
\newcommand{\card}{{\operatorname{Card}}}

\newcommand{\cen}{\operatorname{cen}}

\newcommand{\covol}{\operatorname{Covol}}

\newcommand{\diam}{{\operatorname{diam}}}

\newcommand{\haarhheis}{\operatorname{Haar}_{\operatorname{\HHeis}_7}}

\newcommand{\HHeis}{\HH\!\operatorname{eis}}

\newcommand{\id}{\operatorname{id}}
\renewcommand{\Im}{{\operatorname{Im}}}

\newcommand{\Vol}{\operatorname{Vol}}
\newcommand{\vol}{\operatorname{vol}}

\newcommand{\hdh}{{\HH}^2_\HH}
\newcommand{\hnh}{{\HH}^n_\HH}

\newcommand{\GL}{\operatorname{GL}}
\newcommand{\Sp}{\operatorname{Sp}}

\newcommand{\PU}{\operatorname{PU}}

\newcommand{\tr}{\operatorname{\tt tr}}
\newcommand{\n}{\operatorname{\tt n}}

\newcounter{fig}

\title{Rigidity, counting and equidistribution \\
of quaternionic Cartan chains}
\author{Jouni Parkkonen \and Fr\'ed\'eric Paulin}

\begin{document}
\bibliographystyle{../alphanum}
\maketitle
\begin{abstract}
In this paper, we prove an analog of Cartan's theorem, saying that the
chain-preserving transformations of the boundary of the quaternionic
hyperbolic spaces are projective transformations. We give a counting
and equidistribution result for the orbits of arithmetic chains in the
quaternionic Heisenberg group.
\end{abstract}

\section{Introduction}
\label{sect:intro}

The sphere at infinity $\partial_\infty X$ of a negatively curved
symmetric space $X$ carries many rich structures, from the
geometric, analytic and arithmetic points of view. When the sectional
curvature is not constant, the possibilities are particularly rich, for
instance with the Carnot-Carathéodory, sub-Riemannian or (hyper) CR
structures (see for instance \cite{Mostow73,Gromov96,Goldman99,
  Biquard01, KamNay13}), leading to strong rigidity properties, as
Pansu's rigidity theorem for quasi-isometries \cite{Pansu89}.
Arithmetic subgroups of the isometry group of $X$ endow the sphere at
infinity of $X$ with arithmetic structures, and problems of
equidistribution of rational points or subvarieties in
$\partial_\infty X$, as well as in other homogeneous manifolds, have
been intensively studied (see for instance \cite{Duke03, GorMau05,
  BenOh12, EllMicVen13, BenQui13a, Kim15, BroParPau19, ParPau20a} and
many others).

In this paper, we study the quaternionic hyperbolic spaces $X$, whose
extreme rigidity is exemplified by the Margulis-Gromov-Schoen theorem
in \cite{GroSch92}, proving, contrarily to the real or complex case,
the arithmeticity of lattices in the isometry group of $X$. As
announced in \cite{ParPau20a}, we prove a von Staudt-Cartan type of
rigidity result for the family of all $3$-sphere chains in the sphere
at infinity of $X$, and, analogously to the complex hyperbolic case
treated in \cite{ParPau17MA}, an effective equidistribution result for
the arithmetic chains in orbits of arithmetic groups built using
maximal orders in rational quaternion algebras.

\medskip
More precisely, let $\HH$ be Hamilton's quaternion algebra over $\RR$,
with $x\mapsto \overline{x}$ its conjugation, $\n: x\mapsto
x\overline{x}$ its reduced norm, $\tr: x\mapsto x+\overline{x}$ its
reduced trace.  Let $q$ be the quaternionic Hermitian form on the
right vector space $\HH^{3}$ over $\HH$ defined by
$$
q(z_0,z_1, z_2)=-\tr(\,\overline{z_0}\,z_2) + \n(z_1)\;,
$$ 
and $\PU_q$ its projective unitary group. It is the isometry group
of the quaternionic hyperbolic plane $\hdh$, realised as the negative
cone of $q$ in the right projective plane $\PP^2_{\rm r} (\HH)$, and
normalised to have maximal sectional curvature $-1$.  See Section
\ref{sect:quathypheis} for a more complete description.

The boundary at infinity $\partial_\infty\hdh$ of $\hdh$ is the
isotropic cone of $q$ in $\PP^2_{\rm r} (\HH)$, and the intersections
with $\partial_\infty\hdh$ of the quaternionic projective lines
meeting $\hdh$ are called {\it chains}. We study them, giving their
elementary properties and complete geometric descriptions in Section
\ref{sect:chain}. Our first result is similar to Cartan's theorem (see
\cite{Cartan32,Goldman99}) in the complex hyperbolic case. See Theorem
\ref{theo:Cartan} for a version in any dimension.

\btheo\label{theo:cartanintro}
A chain-preserving transformation from the boundary at infinity of the
quaternionic hyperbolic plane to itself is a projective unitary
transformation.
\etheo

The boundary at infinity $\partial_\infty\hdh$ of $\hdh$, with the
point $\infty=[1:0:0]$ removed, identifies by the map $(w_0,w)\mapsto
[w_0:w:1]$ with the quaternionic Heisenberg group
$$
\HHeis_{7}= \{(w_0,w)\in\HH\times\HH\;:\; \tr\;w_0=\n(w)\}\,,
$$ 
with group law
\begin{equation}\label{eq:heislaw}
(w_0,w)(w'_0,w')= (w_0+w'_0+\,\overline{w}\,w',w+w')\;.
\end{equation} 
We endow the metabelian simply connected real Lie group $\HHeis_{7}$
with its Cygan distance $d_{\rm Cyg}$, which is the unique
left-invariant distance such that $d_{\rm Cyg}((w_0,w),(0,0))=
(4\n(w_0))^{\frac{1}{4}}$.  The chains $C$ contained in $\HHeis_7$ are
ellipsoids, and have a natural center $\cen(C)$ and radius (see
Section \ref{sect:chain}).

Let $A$ be a definite ($A\otimes_\QQ\RR=\HH$) quaternion algebra over
$\QQ$, with discriminant $D_A$. Let $\OOO$ be a maximal order in $A$.
We refer for instance to \cite{Vigneras80} for background on
quaternion algebras and orders.  The group $\PU_q(\OOO)$ of elements
of $\PU_q$ represented by matrices with coefficients in $\OOO$ is a
(necessarily arithmetic) lattice in $\PU_q$.  A chain $C_0$ is said to
be {\it arithmetic} over $\OOO$ if the orbit of some point of $C_0$
under the stabiliser of $C_0$ in $\PU_q(\OOO)$ is dense in $C_0$. The
stabiliser $\PU_q(\OOO)_\infty$ of $[1:0:0]$ in $\PU_q(\OOO)$
preserves the diameters of the chains for $d_{\rm Cyg}$. The following
result (see Theorem \ref{theo:countchain} for an explicit and more
general version) is an asymptotic counting result of the arithmetic
chains in an orbit under the arithmetic group $\PU_q(\OOO)$ when their
Cygan diameter tends to $0$.

\btheo\label{theo:countingintro} Let $C_0$ be an arithmetic chain
in $\partial_\infty\hdh$.  There exists a constant $\kappa>0$ and an
explicit constant $c>0$ such that, as $\epsilon\ra 0$, the number of
chains modulo $\PU_q(\OOO)_\infty$ in the $\PU_q(\OOO)$-orbit of
$C_0$, with Cygan diameter at least $\epsilon$, is equal to
$c\;\epsilon^{-10}(1+\bigO(\epsilon^{\kappa}))$.  
\etheo

An arithmetic chain $C_0$ bounds in $\hdh$ a homothetic copy of the
real hyperbolic space of dimension $4$ normalised to have sectional
curvature $-1$. We denote by $\covol(C_0)$ the volume of the quotient
of this real hyperbolic space by the stabiliser $\PU_q(\OOO)_{C_0}$ of
$C_0$ in $\PU_q(\OOO)$, and by $m_0$ the order of the pointwise
stabiliser of this real hyperbolic space in $\PU_q(\OOO)$. We endow
the real Lie group $\HHeis_{7}$ with its Haar measure $\haarhheis$
normalised in such a way that the total mass of the induced measure on
the quotient of $\HHeis_{7}$ by its (uniform) lattice $\HHeis_7\cap
(\OOO\times\OOO)$ is $\frac{D_A^2}{4}$ (see for instance
\cite[Lem.~8.4]{ParPau20a} for an explanation of this normalisation).
Let $m_A=72$ if $D_A$ is even, and $m_A=1$ otherwise.
Finally, we denote by $\Delta_x$ the unit Dirac mass at any point
$x$. The following result proves that the centers of the arithmetic
chains in an orbit under the arithmetic group $\PU_q(\OOO)$
equidistribute in the quaternionic Heisenberg group.

\btheo\label{theo:equidistribintro} For the weak-star convergence of
measures on $\HHeis_7$, we have
$$
\frac{m_0\,m_A\,\pi^6\,\;\prod_{p|D_A}(p-1)(p^2+1)(p^3-1)}
     {25515\;2^{24}\covol(C_0)} \;\epsilon^{10}\!
\sum_{\scriptsize\begin{array}{c}[g]\in  \PU_q(\OOO)/\PU_q(\OOO)_{C_0}\\
\epsilon\leq \diam_{d_{\rm Cyg}}(gC_0)<\infty\end{array}}\!\!
\Delta_{\cen(gC_0)}\;\weakstar\;\haarhheis.
$$
\etheo

We refer to Section \ref{sect:countequidistrib} for a version with
congruences and error terms, and a more developped study of explicit
examples of arithmetic chains.

\medskip
\noindent{\small {\it Acknowledgements: } The authors thank the snowy
  artic conditions in \"Ak\"aslompolo in January 2020 which have
  provided an exceptional working environment. This research was
  supported by CNRS IEA BARP. The second author thanks the Laboratoire
  de Mathématiques Jean Leray at the Université de Nantes where this
  paper was completed.}

\section{Quaternionic hyperbolic spaces and Heisenberg groups}
\label{sect:quathypheis}

In this section, we briefly recall some background on the quaternionic
hyperbolic spaces and quaternionic Heisenberg group, as mostly
contained in \cite[\S 2 and \S 6]{ParPau20a}, see also
\cite{KimPar03,Philippe16} (with different choices of quaternionic
Hermitian form and normalisation of the curvature).

Let $\HH$ be Hamilton's quaternion algebra over $\RR$, with $x\mapsto
\overline{x}$ its conjugation, $\n: x\mapsto x\overline{x}$ its
reduced norm, $\tr: x\mapsto x+\overline{x}$ its reduced trace. We
denote by $(1,i,j,k)$ the canonical basis of $\HH$ as a real vector
space, so that $\overline{x_0+x_1\,i+x_2\,j+x_3\,k}=
x_0-x_1\,i-x_2\,j-x_3\,k$. Let
$$
\Im\;\HH=\{x\in\HH\;:\;\tr x=0\}= \RR\,i+\RR\,j+\RR\,k
$$
be the $\RR$-subspace of purely imaginary quaternions of $\HH$. For
all $w=(w_1,\dots,w_{N})$ and $w'=(w'_1,\dots,w'_{N})$ in $\HH^{N}$,
we denote by $\overline{w}\cdot w'= \sum_{p=1}^{N}
\overline{w_p}\,w'_p$ their standard quaternionic Hermitian product,
and we define $\n(w)= \overline w\cdot w= \sum_{p=1}^{N} \n(w_p)$. We
endow $\HH^{N}$ with the standard Euclidean structure $(w,w')\mapsto
\frac12 \tr(\,\overline{w}\cdot w')$.

We fix $n\in\NN-\{0,1\}$. On the right vector space
$\HH\times\HH^{n-1}\times\HH$ over $\HH$ with coordinates
$(z_0,z,z_n)$, let $q$ be the nondegenerate quaternionic Hermitian
form
\begin{equation}\label{eq:q}
q(z_0,z,z_n)=-\tr(\,\overline{z_0}\,z_n) + \n(z)
\end{equation}
of Witt signature $(1,n)$, and let $\Phi:\HH^{n+1}\times \HH^{n+1}\ra
\HH$, defined by
\begin{equation}\label{eq:Phi}
\Phi:((z_0,z,z_n),(z'_0,z',z'_n))\mapsto
-\overline{z_0}\,z'_n-\overline{z_n}\,z'_0+\,\overline{z}\cdot z'\,,
\end{equation}
be the associated quaternionic sesquilinear form.

The {\it Siegel domain} model of the quaternionic hyperbolic $n$-space
$\hnh$ is
$$ 
\big\{(w_0,w)\in\HH\times\HH^{n-1}\;:\; 
\tr\, w_0 -\n(w)>0\big\}\,,
$$ 
endowed with the Riemannian metric
$$
ds^2_{\,\hnh}=\frac{1}{(\tr\, w_0 -\n(w))^2}
\big(\,\n(dw_0-\overline{dw}\cdot w)+
(\tr\, w_0-\n(w))\;\n(dw)\,\big)\,.
$$
Its boundary at infinity is
$$
\partial_\infty\hnh=\big\{(w_0,w)\in
\HH\times \HH^{n-1} \;:\; \tr\, w_0 -\n(w)=0\big\}\cup\{\infty\}\,.
$$
A {\it quaternionic geodesic line} in $\hnh$ is the image by an
isometry of $\hnh$ of the intersection of $\hnh$ with the quaternionic
line $\HH\times\{0\}$. With our normalisation of the metric, a
quaternionic geodesic line is a totally geodesic submanifold of real
dimension $4$ and constant sectional curvature $-4$.

The closed horoballs in $\hnh$ centred at $\infty\in\partial_\infty
\hnh$ are the subsets
\begin{equation}\label{eq:defhoroinfty}
\H_s=\{(w_0,w)\in\hnh\;:\;\tr w_0-\n(w)\geq s\},
\end{equation}
and the horospheres centred at $\infty$ are their boundaries
$\partial\H_s$, where $s$ ranges in $]0,+\infty[\,$.  Note that, for
every $s\in\;]0,1]$, we have
\begin{equation}\label{eq:distentrhorob}
d(\partial\H_1,\partial\H_s)=-\frac{\ln s}{2}\;.
\end{equation}

The Siegel domain $\hnh$ embeds in the right quaternionic projective
$n$-space $\PP^n_{\rm r}(\HH)$ by the map (using homogeneous
coordinates)
$$
(w_0,w)\mapsto [w_0:w:1]\;.
$$ 
By this map, we identify $\hnh$ with its image, which when endowed
with the isometric Riemannian metric, is called the {\it projective
  model} of $\hnh$. Note that this image is the {\it negative cone} of
the quaternionic Hermitian form $q$ defined in Equation \eqref{eq:q}~:
we have $\hnh=\big\{[z_0:z:z_n]\in \PP^n_{\rm r}(\HH) \;:\;
q(z_0,z,z_n) <0\big\}$.  This embedding extends continuously to the
boundary at infinity, by mapping the point $(w_0,w)\in \partial_\infty
\hnh -\{\infty\}$ to $[w_0:w:1]$ and $\infty$ to $[1:0:0]$, so that
the image of $\partial_\infty\hnh$ is the {\it isotropic cone} of $q$:
we have $\partial_\infty\hnh=\big\{[z_0:z:z_n]\in\PP^n_{\rm r}(\HH)
\;:\; q(z_0,z,z_n)=0 \big\}$. A projective point $[z_0:z:z_n]\in
\PP^n_{\rm r}(\HH)$ is {\it positive} if $q(z_0,z,z_n)>0$.

For every $N\in\NN$, let $I_N$ be the identity $N\times N$ matrix. Let
$$
J=\begin{pmatrix} 0 & 0 & -1\\0 & I_{n-1} & 0\\-1 & 0 & 0
\end{pmatrix}\;.
$$
The {\it conjugate-transpose} matrix a quaternionic matrix $X=
(x_{p,p'})_{1\leq p\leq r,\, 1\leq p'\leq s} \in \M_{r,s}(\HH)$ is $X^*=
(x^*_{p,p'}= \overline{x_{p',p}}\,)_{1\leq p\leq s,\,1\leq p'\leq r}
\in \M_{s,r}(\HH)$.  Let
$$
\operatorname{U}_q=\{g\in \GL_{n+1}(\HH)\;:\; q\circ g=q\}
=\{g\in \GL_{n+1}(\HH)\;:\; g^*J\,g=J\}
$$ 
be the {\it unitary group} of $q$.  Its left linear action on $\HH^{n+1}$
induces a projective action on $\PP^n_{\rm r}(\HH)$ with kernel its
center, which is reduced to $\{\pm I_{n+1}\}$.  The {\it projective
  unitary group}
$$
\PU_q=\operatorname{U}_q/\{\pm I_{n+1}\}
$$ 
of $q$ acts faithfully on $\PP^n_{\rm r}(\HH)$, preserving $\hnh$,
and its restriction to $\hnh$ is the full isometry group of $\hnh$.

A matrix
$$
X=\begin{pmatrix} a & \gamma^* & b\\ \alpha & A & \beta 
\\ c & \delta^* & d\end{pmatrix}\in \GL_{n+1}(\HH)\,,
$$
with $a,b,c,d\in\HH$, $\alpha,\beta,\gamma,\delta \in\HH^{n-1}$
(identified with their column matrices in $\M_{n-1,1}(\HH)$) and
$A\in\M_{n-1,n-1}(\HH)$, belongs to $\operatorname{U}_{q}$ if and only
if
\begin{equation}\label{eq:equationsUq}
\begin{cases}
\hfill \overline{c}\,a-\alpha^*\alpha+\overline{a}\,c   &= 0 \\ 
\hfill  \overline{d}\,b-\beta^*\beta+\overline{b}\,d   &= 0 \\ 
   -\delta\,\ga^*+A^* A - \ga\,\delta^*  \!\!\!&=I_{n-1} \\
\hfill  \overline{d}\,a-\beta^*\alpha+\overline{b}\,c &= 1\\
\hfill  \delta \,a-A^*\alpha+\ga\,c &=0\\
\hfill  \delta\,b-A^*\beta+\ga\,d  &=0 \;.
\end{cases}
\end{equation}
With $\Sp(n-1)=\{g\in \GL_{n+1}(\HH): g^*g=I_{n-1}\}$, an easy
computation shows that the block upper triangular subgroup of
$\operatorname{U}_q$ is
$$
\operatorname{B}_q=\Bigg\{\begin{pmatrix} 
\mu r & \zeta^* & \frac{1}{2r}(\n(\zeta)+u)\mu \\
0 & U & \frac{1}{r}\;U\,\zeta\;\mu\\
0 & 0 & \frac{\mu}{r}
\end{pmatrix}\;:\; 
\begin{array}{c}\zeta\in \HH^{n-1},\; u\in\Im\;\HH,\\
U\in \Sp(n-1), \mu\in \Sp(1), r>0\end{array}\Bigg\}\;.
$$ 
Its image $\operatorname{PB}_q=\operatorname{B}_q/\{\pm I_{n+1}\}$
in $\PU_q$ is equal to the stabiliser of $\infty$ in $\PU_q$.

\medskip
The {\it quaternionic Heisenberg group} $\HHeis_{4n-1}$ of dimension
$4n-1$ is the real Lie group structure on $\HH^{n-1}\times \Im\;\HH$
with law
$$
(\zeta,u)(\zeta',u')=
(\zeta+\zeta',u+u'+2\,\Im\;\overline{\zeta}\cdot\zeta')
$$ 
and inverses $(\zeta,u)^{-1}=(-\zeta,-u)$. 
It
identifies with $\partial_\infty\hnh-\{\infty\}$ by the map
$(\zeta,u)\mapsto (w_0,w)$ where
\begin{equation}
\label{eq:changecoord}
(w_0,w)  =\Big(\;\frac{\n(\zeta)+u}{2},\;\zeta\,\Big)\;\;\;
{\rm hence}\;\;\;
(\zeta,u)  =(w,\;2\,\Im\; w_0)\,,
\end{equation}
and with a subgroup of $\operatorname{PB}_q \subset\PU_q$ preserving
every horoball $\H_s$ for $s>0$ by the map $(\zeta,u)\mapsto
\pm\begin{pmatrix} 1 & \zeta^* & \frac{\n(\zeta)+u}{2}\\0 & I_{n-1} &
\zeta\\0 & 0 & 1 \end{pmatrix}$. Equation \eqref{eq:changecoord}
allows to recover the definition of $\HHeis_{7}$ given in the
Introduction, for which the inverses are $(w_0,w)^{-1}=
(-w_0+\n(w),-w)$.

For every $(\zeta,u)\in\HHeis_{4n-1}$, the map $(\zeta',u')\mapsto
(\zeta,u) (\zeta',u')$ is the {\em Heisenberg translation} by
$(\zeta,u)$.  For every $\zeta\in\HH^{n-1}$, the Heisenberg
translation by $(\zeta,0)$ is called a {\it horizontal (Heisenberg)
  translation}. For every $u\in\Im\;\HH$, the Heisenberg translation
by $(0,u)$ is called a {\it vertical (Heisenberg) translation}. The
canonical map $\Pi_v:\HHeis_{4n-1} \ra\HH^{n-1}$ defined by $(\zeta,u)
\mapsto \zeta$ is a real Lie group morphism, called the {\it vertical
  projection}, whose kernel is the center of $\HHeis_{4n-1}$. For
every $U\in\Sp(n-1)$, the map $(\zeta,u)\mapsto (U\zeta,u)$ is the
{\em Heisenberg rotation} by $U$. For every $\lambda>0$, the map
$h_\lambda: (\zeta,u)\mapsto (\lambda\zeta,\lambda^2u)$ is the {\em
  Heisenberg dilation} by $\lambda$.

The {\it Cygan distance} $d_{\rm Cyg}$ on $\HHeis_{4n-1}$ is the
unique left-invariant distance on the real Lie group $\HHeis_{4n-1}$
such that
\begin{equation}\label{eq:defidistCyg}
d_{\rm Cyg}((\zeta,u),(0,0))=\big(\,\n(\zeta)^2+ \n(u)\,\big)^{1/4}\,,
\end{equation}
or equivalently $d_{\rm Cyg}((w_0,w),(0,0))= (4\n(w_0))^{\frac{1}{4}}$
by Equation \eqref{eq:changecoord}. We introduce (see
\cite{ParPau11MZ,ParPau17MA} in the complex case) the {\it modified
  Cygan distance} $d''_{\rm Cyg}\,$, as the unique left-invariant map
from $\HHeis_{4n-1}\times\HHeis_{4n-1}$ to $[0,+\infty[$ such that
\begin{equation}\label{eq:defidistCyg''}
d''_{\rm Cyg}((\zeta,u),(0,0))=
\frac{(\n(\zeta)^2+ \n(u))^{1/2}}
{\big(\,(\n(\zeta)^2+ \n(u))^{1/2}+\n(\zeta)\,\big)^{1/2}}\,,
\end{equation}
or equivalently by Equation \eqref{eq:changecoord}
$$
d''_{\rm Cyg}((w_0,w),(0,0))=
\frac{2\n(w_0)^{1/2}}{(2\n(w_0)^{1/2}+\n(w))^{1/2}}\,.
$$
Though not actually a distance, the map $d''_{\rm Cyg}$ is symmetric
and satisfies
$$
\frac{1}{\sqrt{2}}\,d_{\rm Cyg}\leq d''_{\rm Cyg}\leq d_{\rm Cyg}\,.
$$ 
For every nonempty bounded subset $A$ of $\HHeis_{4n-1}$, we define
the {\it diameter} of $A$ for this almost distance as
$$
\diam_{d''_{\rm Cyg}}(A)=\sup_{x,\,y\,\in A} d''_{\rm Cyg} (x,y)\;.
$$
Note that the Cygan distance and the modified Cygan distance are
invariant under Heisenberg translations and rotations, and that for
every $\lambda>0$, the Heisenberg dilation $h_\lambda$ is a homothety
of ratio $\lambda$ for both distances.

\blemm\label{lem:distgeod} 
For every geodesic line $]x,y[$ in $\hnh$ disjoint from the horoball
$\H_1$, the distance in $\hnh$ between $\H_1$ and $]x,y[$ is equal to
$$
d(\H_1,\,]x,y[)=
-\ln \Big(\frac{1}{\sqrt{2}} \,d''_{\rm Cyg}(x,y)\Big)\,.
$$
\elemm

\dem By the invariance under Heisenberg translations of $\H_1$, of the
distance in $\hnh$ and of the modified Cygan distance, we may assume
that $x=(w_0,w)\in \partial_\infty\hnh-\{\infty,(0,0)\}$ and
$y=(0,0)\in \partial_\infty\hnh -\{\infty\}$. By
\cite[Lem.~6.4]{ParPau20a}, the geodesic line from $(w_0,w)$ to
$(0,0)$ is, up to translation at the source, the map
$$
\ga_{w_0,w}:t\mapsto \big( w_0(1+e^{2t}w_0)^{-1},\,w(1+e^{2t}w_0)^{-1}\,)\,.
$$ 
The point $\gamma_{w_0,w}(t)$ belongs to the horosphere $\H_{s(t)}$,
where, since $\tr w_0=\n(w)$,
$$
s(t)= \tr( w_0(1+e^{2t}w_0)^{-1})-\n(w(1+e^{2t}w_0)^{-1})
= \frac{2\,e^{2t}\n(w_0)}{\n(1+e^{2t}w_0)}.
$$ 
Let $r=\n(w_0)^{1/2}$ be the norm of the vector $w_0$ and $\theta$
the angle between the vectors $1$ and $w_0$ in the Euclidean space
$\HH$. Then the map 
$$
t\mapsto s(t)=\frac{2\,e^{2t}r^2}{e^{4t}r^2+2r e^{2t}\cos\theta+1}
$$ 
reaches it maximum at $e^{2t}=\frac{1}{r}$. Since $\tr w_0=\n(w)$,
the value of this maximum is
$$
s_{\max}=\frac{2\n(w_0)^{1/2}}{2+\tr(w_0\n(w_0)^{-1/2})}=
\frac{2\n(w_0)}{2\n(w_0)^{1/2}+\n(w)}=
\frac{1}{2}\,d''_{\rm Cyg} ((w_0,w),(0,0))^2\,.
$$
The result then follows from Equation \eqref{eq:distentrhorob}.
 \cqfd

\section{Chains}
\label{sect:chain}

In this section, we define the quaternionic Cartan chains and give
their elementary geometric properties, see also \cite{Shi16}.  In the
complex case, the notion of chain is attributed to von Staudt by
\cite{Cartan32}.  The exposition follows the one of \cite{Goldman99}
in the complex case. We fix $m\in\{1,\dots, n-1\}$.

\subsection{A vocabulary of chains}
\label{subsec:vocabulary}

An {\it $m$-chain} $C$ in $\partial_\infty\hnh$ is the intersection
with $\partial_\infty\hnh$ of a quaternionic projective space $L_C$ of
dimension $m$ meeting $\hnh$. Note that $C$ determines $L_C$ and
conversely. A {\it chain} is a $1$-chain, and a {\it hyperchain} is an
$(n-1)$-chain. A $m$-chain is {\it vertical} if it contains $\infty
=[1:0:0]$, and {\it finite} otherwise.

If $P=[z_0:z:z_{n}]\in\PP^{n}_{\rm r}(\HH)$, let 
$$
P^\perp= \big\{[z'_0:z':z'_{n}] \in \PP^{n}_{\rm r}(\HH): 
\Phi((z_0,z,z_{n}), (z'_0,z',z'_{n}))=0\big\}
$$
be the orthogonal quaternionic projective subspace of $P$. The map
$P\mapsto P^\perp$, from the set of positive projective points to the
set of quaternionic projective hyperplanes in $\PP^{n}_{\rm r}(\HH)$
meeting $\hnh$, is a $\PU_q$-equivariant bijection.  Therefore, the
map
$$
P\mapsto C_P=P^\perp\cap \partial_\infty \hnh
$$
is a $\PU_q$-equivariant bijection from the set of positive projective
points to the set of hyperchains. The point $P$ is called the {\it
  polar point} of the hyperchain $C_P$, or of the quaternionic
projective hyperplane $P^\perp$. If $P=[z_0:z:z_{n}]$, we have
\begin{equation}\label{eq:eqchain}
C_P\cap(\partial_\infty\hnh-\{\infty\})=\big\{[w_0:w:1]\; : 
-\big(\frac{\n(w)}{2}-\Im \,w_0\big)z_n+
\overline{w}\cdot z-z_0=0\big\}\;.
\end{equation}
This hyperchain $C_P$ is hence vertical if and only if $z_n=0$, in
which case $C_P\cap(\partial_\infty\hnh-\{\infty\})$ is the preimage
by the vertical projection $\Pi_v:\HHeis_{4n-1}\ra \HH^{n-1}$ of the
quaternionic affine hyperplane of $\HH^{n-1}$ with equation
$\overline{z}\cdot w=\overline{z_0}$ in the unknown $w$. Similarly a
vertical chain is the preimage of a point of $\HH^{n-1}$ by the
vertical projection $\Pi_v$.

When $C=C_P$ is a finite hyperchain, that is, when $z_n\neq 0$, then
$C$ is a codimension $4$ ellipsoid in the Euclidean space
$\HH^{n-1}\times\Im\,\HH$, whose vertical projection is the Euclidean
sphere in $\HH^{n-1}$ with real codimension $1$ and equation
$\n(w)-\tr (\overline{w}\cdot zz_n^{-1})+\tr(z_0z_n^{-1})=0$ in the
unknown $w$, with center $zz_n^{-1}$ and radius
$$
R_C=\frac{q(z_0,z,z_n)^{1/2}} {\n(z_n)^{1/2}}\;.
$$
This radius $R_C$ of the Euclidean sphere $\Pi_v(C)$ is called the
{\it radius} of the finite hyperchain $C$. The map $\Pi_v|_{C}$ from
$C$ to $\Pi_v(C)$ is a homeomorphism. When $z=0$ and
$z_0z_n^{-1}\in\RR$, the hyperchain $C=C_P$ is contained in the
horizontal subspace $\{(\zeta,u)\in\HH^{n-1}\times \Im\,\HH: u=0\}$ of
$\HHeis_{4n-1}$, by Equation \eqref{eq:changecoord}.

Similarly, a finite chain is a $3$-dimensional ellipsoid in the
Euclidean space $\HH^{n-1} \times\Im\,\HH$, whose vertical projection
is a Euclidean $3$-sphere in $\HH^{n-1}$. In particular, any chain is
homeomorphic to the $3$-sphere $\SSS^{3}$.

\subsection{Transitivity properties of $\PU_q$ on chains}
\label{subsec:transitivity}

Through any two distinct projective points belonging to $\partial_\infty
\hnh$ passes one and only one quaternionic projective line, and this
projective line meets $\hnh$. Hence through two distinct points of
$\partial_\infty\hnh$ passes one and only one chain.
By Witt's theorem, the group $\PU_q$ acts transitively on the set of
quaternionic projective spaces $L$ of dimension $m$ meeting $\hnh$,
hence it acts transitively on the set of $m$-chains. 

Note that two $m$-chains having the same vertical projection differ by
a vertical Heisenberg translation, that the group generated by
Heisenberg translations and rotations acts transitively on the set of
vertical $m$-chains, and that $\operatorname{PB}_q$ (that contains the
Heisenberg dilations, rotations and translations) acts transitively on
the set of finite $m$-chains.

The next result gives the topological structure of a family of chains,
called a {\it fan} in the complex hyperbolic case (see for instance
\cite[page 131]{Goldman99}).

\bprop The union $F$ of all chains containing a given point
$P\in\partial_\infty \hnh$ and passing through an $m$-chain $C$ of
$\partial_\infty \hnh$ not containing $P$ is homeomorphic to the
topological quotient space $(\SSS^3\times\SSS^{4m-1})/\!\!\sim$ where
$\sim$ is the equivalence relation generated by $(x_0,x)\sim (x_0,y)$
for all $x,y\in\SSS^{4m-1}$, where $x_0$ is any fixed point in
$\SSS^3$.
\eprop

\dem By the transitivity properties of $\PU_q$, we may assume that
$P=\infty$. Hence $C$ is a finite chain, and by the transitivity
properties of the Heisenberg translations, we may assume that $C$ is a
Euclidean sphere of dimension $4m-1$ contained in the horizontal space
$\{(\zeta,u)\in\HH^{n-1}\times \Im\,\HH: u=0\}$. Thus $F=
\bigcup_{(\zeta,u)\in C}\Pi_v^{-1}(\zeta,u)$ is clearly homeomorphic
to the above quotient of $\SSS^3\times\SSS^{4m-1}$.
\cqfd

\subsection{Reflexions on chains}
\label{subsec:reflexions}

The chains are fixed point sets at infinity of natural isometries of
$\hnh$, that we now describe.

If $L$ is a proper quaternionic projective subspace of $\PP^n_{\rm r}(\HH)$,
there exists a unique involution $\iota_L$ in $\PU_q$ with fixed point
set $L$, called the {\it reflexion} on $L$. Note that the set of fixed
point of $\iota_L$ in $\partial_\infty\hnh$ is the $m$-chain
$L\cap\partial_\infty\hnh$, where $m$ is the quaternionic dimension of
$L$, assuming that $m\neq 0$.

For instance, $C=\big\{[z_0:z_1:\dots:z_{n}] \in\partial _\infty\hnh:
z_{m}=0, \dots, z_{n-1}=0\big\}\cup\{\infty\}$ is a vertical
$m$-chain, called the {\it standard vertical $m$-chain} and the
reflexion $\iota_{L_C}$ is the map
$$
[z_0:z_1:\dots :z_{n}]\mapsto [z_0:z_1:\dots :z_{m-1}:
  -z_m: \dots :-z_{n-1}:z_{n}]\,.
$$
The vertical $m$-chains are the images of the standard vertical
$m$-chain by the Heisenberg translations and Heisenberg rotations:
they are the
$$
(E\times\Im\;\HH)\cup\{\infty\}
$$
where $E$ is a quaternionic affine subspace of $\HH^{n-1}$ with
dimension $m-1$ (hence a point when $m=1$).

\blemm \label{lem:reflexcommut} Let $L$ and $L'$ be quaternionic
projective subspaces of $\PP^n_{\rm r}(\HH)$ meeting $\hnh$ such that
one is not contained in the other, whose sum of dimensions is $n$.
The following assertions are equivalent.
\begin{enumerate}
\item[(1)] The reflexions $\iota_L$ and $\iota_{L'}$ commute.
\item[(2)] The reflexion $\iota_L$ preserves $L'$.
\item[(3)] The reflexion $\iota_{L'}$ preserves $L$.
\item[(4)] We have $(\iota_L\circ \iota_{L'})^2=\id$.
\item[(5)] The totally geodesic subspaces $L\cap\hnh$ and $L'\cap\hnh$
  intersect perpendicularly in the Riemannian manifold $\hnh$.
\item[(6)] The  subspace $L\cap\hnh$ is a fiber
  of the orthogonal projection on $L'\cap\hnh$ in $\hnh$.
\item[(7)] The subspace $L'\cap\hnh$ is a fiber
  of the orthogonal projection on $L\cap\hnh$ in $\hnh$.
\end{enumerate}
\elemm

\dem The proof is similar to the one of \cite[Lem.~4.3.1]{Goldman99}
in the complex hyperbolic case. Note that $L\cap\hnh$, being the set
of fixed points of the isometry $\iota_L$ of the negatively curved
Riemannian manifold $\hnh$, is indeed totally geodesic.

Two involution commutes if and only if their composition is an
involution or the identity, hence Assertions (1) and (4) are
equivalent.  Since the centralizer of a projective transformation
preserves its fixed point set, Assertion (1) implies Assertions (2)
and (3).  If Assertion (2) is satisfied, then
$\iota_L\circ\iota_{L'}\circ\iota_{L}^{\;-1}= \iota_{\iota_{L}(L')} =
\iota_{L'}$, so that Assertion (1) is satisfied. Similarly, Assertion
(2) implies Assertion (1).  Finally, the totally geodesic subspaces
$L'\cap\hnh$ and $L\cap\hnh$ in $\hnh$

$\bullet$~ either have disjoint closures in $\hnh\cup\partial_\infty
\hnh$,

$\bullet$~ or are disjoint and have closures meeting in
$\partial_\infty\hnh$,

$\bullet$~ or meet in $\hnh$.

\noindent In the first two cases, the composition $\iota_L\circ
\iota_{L'}$ has infinite order, and in the last case, $\iota_L\circ
\iota_{L'}$ can be an involution if and only if $L'\cap\hnh$ and
$L\cap\hnh$ are perpendicular.
\cqfd

\medskip
An $m$-chain $C$ and an $(n-m)$-chain $C'$ are {\it orthogonal} if
neither of the corresponding quaternionic projective subspaces $L_C$
and $L_{C'}$ contains the other and if they satisfy one of the
equivalent assertions of Lemma \ref{lem:reflexcommut}. For instance,
the hyperchains orthogonal to the standard vertical chain
$(\{0\}\times \Im\;\HH)\cup\{\infty\}$ are exactly the Euclidean
spheres centered at $(0,u_0)$ in the horizontal subspace $\{(\zeta,u)
\in\HH^{n-1}\times \Im\,\HH: u=u_0\}$ of $\HHeis_{4n-1}$, for some
$u_0$ in $\Im\;\HH$.

\subsection{Description of the center and radius of chains}
\label{subsec:centers}

We now define and study the centers of chains, whose equidistribution
we will prove in Section \ref{sect:countequidistrib}. 

The {\it center} of an $m$-chain $C$ is $\cen(C)=\iota_{L_C}(\infty)$.
In particular, $\cen(C)=\infty$ if and only if $C$ is vertical.  For
every element $\ga\in \operatorname{PB}_q$ (which fixes $\infty$), the
reflexion on the finite $m$-chain $\ga C$ is $\ga\iota_{L_C}\ga^{-1}$,
so that the center of $\ga C$ is
\begin{equation}\label{eq:equivcenter}
\cen(\ga C)=\ga \cen (C)\,.
\end{equation}

When $P_0=[-\frac{1}{2}:0:1]$, the hyperchain $C_{P_0}$ with polar
point $P_0$ is, by Equation \eqref{eq:eqchain}, the sphere centered at
$(0,0)$ with radius $1$ in the horizontal codimension $3$ Euclidean
subspace $\{(\zeta,u)\in\HH^{n-1}\times \Im\,\HH: u=0\}$ in
$\HHeis_{4n-1}$. The reflexion on $L=L_{C_{P_0}}$ is the involutive map
$\iota_L:(w_0,w) \mapsto (\frac{1}{4}\,w_0^{\;-1},
\frac{1}{2}\,ww_0^{\;-1})$, induced by $\pm\begin{pmatrix} 0 & 0 &
1/2 \\ 0 & I_{n-1} & 0\\ 2 & 0 &
0 \end{pmatrix}\in\PU_q$. Thus,
$\cen(C_{P_0})=\iota_L(\infty)=(0,0)$.

Let $P=[z_0:z:z_n]$ be a positive projective point with $z_n\neq
0$. An easy computation shows that the Heisenberg translation $\ga$ by
$$
\big[\frac{\n(z)}{2\n(z_n)}-\Im(z_0z_n^{\;-1}): -zz_n^{\;-1}:1\big]
$$
maps $P$ to $[-\frac{R^2}{2}:0:1]$ where $R=R_{C_P}=
\frac{q(z_0:z:z_n)^{1/2}} {\n(z)^{1/2}}$ is the radius of the finite
hyperchain $C_P$, and the Heisenberg dilation
$$
h_{R}:(w_0,w)\mapsto (R^2w_0,Rw)
$$
maps $P_0$ to $[-\frac{R^2}{2}:0:1]$. Hence the center of the finite
hyperchain $C_P$ with polar point $P$ is, by Equation
\eqref{eq:equivcenter}, equal to
$$
\cen(C_P)=\ga^{-1}h_{R}\cen(C_{P_0})=\ga^{-1}(0,0)=
\big[\frac{2\,\Im(z_0\,\overline{z_n})+\n(z)}{2\n(z_n)}:
  zz_n^{\;-1}:1\big]\,,
$$
or $\cen(C_P)=\big(zz_n^{\;-1},2\,\Im(z_0\,z_n^{-1})\big)$ in the
$(\zeta,u)$-coordinates of $\HHeis_{4n-1}$ by Equation
\eqref{eq:changecoord}. Thus, by Equation \eqref{eq:eqchain}, if $C$
is a finite hyperchain in $\HHeis_{4n-1}$ with center $(\zeta_0,u_0)$
and radius $r_0$, then
$$
C=\{(\zeta,u)\in\HH^{n-1}\times\Im\;\HH\;: \;\n(\zeta-\zeta_0)=r_0^{\;2}
\;\;{\rm and}\;\; u=u_0+2\,\Im(\,\overline{\zeta_0}\,\zeta)\}\,.
$$
In particular, a finite hyperchain is uniquely determined by its
center and its radius, and the hyperchains contained in the horizontal
Euclidean space $\{(\zeta,u)\in\HH^{n-1}\times\Im\;\HH:u=0\}$ are
exactly the Euclidean spheres centered at $(0,0)$.

\subsection{A von Staudt-Cartan rigidity theorem}
\label{subsect:Cartan}

The following theorem shows that the chain-preserving transformations
of the boundary of the quaternionic hyperbolic spaces are projective
transformations.  This is a quaternionic version of the result of
Cartan in the complex case (see for instance
\cite[Theo.~4.3.12]{Goldman99}), close to von Staudt's fundamental
theorem of real projective geometry.

\btheo\label{theo:Cartan} A bijection $f$ from $\partial_\infty\hnh$
to itself, mapping chains to chains, is (the restriction to
$\partial_\infty\hnh$ of) an element of $\PU_q$.  
\etheo

\dem Up to composing by an element of $\PU_q$, we may assume that $f$
fixes $\infty=[1:0:0]$. Hence $f$ preserves the set of vertical
chains, which are the ones containing $\infty$. The set of vertical
chains identifies with the horizontal space $\HH^{n-1}$ of the
quaternionic Heisenberg group by the vertical projection $\Pi_v$,
which sends a vertical chain $C$ to the unique point of $\HH^{n-1}$
whose preimage by $\Pi_v$ is $C$. Hence $f$ induces a bijection
$\overline{f}$ from $\HH^{n-1}$ to itself, which sends the vertical
projections of the finite chains to the vertical projections of the
finite chains.

The vertical projections of the finite chains are exactly all the
Euclidean $3$-spheres in $\HH^{n-1}$.  Given two distinct points $x,y$
in $\HH^{n-1}$, the complement of the union of all the Euclidean
$3$-spheres containing $x$ and $y$ is the real affine line containing
$x$ and $y$, with $x$ and $y$ removed. Hence $\overline{f}$ is a
bijection of $\HH^{n-1}$ sending real affine lines to real affine
lines. By the fundamental theorem of real affine geometry, this map is
an affine transformation of $\HH^{n-1}$. Since the affine
transformations of $\HH^{n-1}$ are vertical projections of elements of
the stabiliser $\operatorname{PB}_q$ of $\infty$ in $\PU_q$, up to
composing $f$ by an element of $\operatorname{PB}_q$, we may assume
that $\overline{f}$ is the identity map of $\HH^{n-1}$, and also that
$f(0)=0$.

Let $x\in\partial_\infty\hnh-\{\infty\}$, and let us prove that
$f(x)=x$.  First assume that $\Pi_v(x)\neq 0$. Then the unique chain
$C_x$ passing through $0$ and $x$ is a finite chain, and the vertical
projections of $C_x$ and $f(C_x)$ coincide, since $\overline{f}=\id$.
By the uniqueness of a chain with given vertical projection up to a
vertical translation, since $f(0)=0$, we have $f(C_x)= C_x$. But if
$f(x)\neq x$, then since $f(x)$ and $x$ have the same vertical
projections, the chains $C_x$ and $f(C_x)$ through $0$ would be
different. Hence $f(x)=x$. This is in particular true for all fixed
$x=x_0\neq 0$ in the horizontal space $\HH^{n-1}\times\{0\}$.
Replacing $0$ by such an $x_0$ in the above argument allows to prove
that $f(x)=x$ when $\Pi_v(x)= 0$.  \cqfd

\medskip

A similar proof shows that an injective map $f$ from $\partial_\infty
\hnh$ to itself, such that any three points belong to a same chain if
and only if their images by $f$ belong to a same chain, is the
restriction of an element of $\PU_q$.

\subsection{Relation with the hyper CR structure}
\label{subsect:CR}

In this subsection, we give a characterisation of the chains in terms
of the natural hyper CR structure on $\partial_\infty\hnh$. We refer
for example to \cite{Besse87} and \cite{KamNay13} for background on
hyperkähler manifolds and hyper CR manifolds, respectively.

We endow the manifold $\PP^n_{\rm r}(\HH)$ with its natural
hyperkähler structure, and we denote by $(\II,\JJ,\KK)$ the
corresponding triple of almost complex structures. The boundary at
infinity $W=\partial_\infty\hnh$ is a smooth real hypersurface in the
real manifold $\PP^n_{\rm r}(\HH)$ of real dimension $4n$, and $E=
TW\cap \,\II\, TW \cap\, \JJ\, TW \cap \,\KK\, TW $ is a real
codimension $3$ subbundle of the real tangent bundle $T\PP^n_{\rm r}
(\HH)$, invariant under $\PU_q$, defining a hyper CR structure on $W$.
When $x$ is the point $(0,0)$ in the $(\zeta,u)$-coordinates of
$\HHeis_{4n-1}=\partial_\infty\hnh -\{\infty\}$, then, identifying
$\HH^{n-1}\times\HH$ with its real tangent space, the fiber $E_x$ of
$E$ over $x$ is the horizontal subspace $\{(\zeta,u)\in\HH^{n-1}\times
\HH:u=0\}$.

A {\it calibration} of $E$ is a $1$-form $w$ with values in $\Im
\;\HH$ such that $E=\ker w$. Its {\it Levi form} is $dw$. For
instance, in the $(\zeta,u)$-coordinates of $\HHeis_{4n-1}$,
$$
\omega=du -2\;\Im\;(\,\overline{\zeta}\cdot d\zeta)
$$
is a calibration of $E$ (when restricted to $\partial_\infty\hnh
-\{\infty\}$). An easy computation shows that this calibration is
invariant under Heisenberg translations and rotations: For every such
transformation $\ga$, we have $\ga^*\omega=\omega$. The fact that
$\omega$ is indeed a calibration follows by invariance since $\ker du
=\{(\zeta,u)\in\HH^{n-1}\times\HH:u=0\}$. This calibration $\omega$ is
scaled by the Heisenberg dilations as follows : for every $\lambda>0$,
we have $(h_\lambda)^*\omega=\lambda^2\,\omega$.

In the following result, we denote by $v=v_1\,i+v_2\,j+v_3\,k$ the
standard coordinate in $\Im\;\HH$. We denote by
$\omega_1,\omega_2,\omega_3$ the standard coordinates of the
calibration $\omega$, so that
$$
\omega=\omega_1 \,i +\omega_2\, j +\omega_3 \,k\;.
$$

Given a chain $C$ in $\partial\hnh$, let $\mu=\mu_C$ be the (Borel
positive) measure on $\HHeis_{4n-1}$ with support $C\cap\HHeis_{4n-1}$
associated with the volume form $\omega_1\wedge \omega_2\wedge
\omega_3$ on $C$. For instance, if
$C=\{(\zeta,u)\in\HH^{n-1}\times\,\Im\;\HH:\zeta=0\}\cup \{\infty\}$
is the standard vertical chain, then $\omega|_C = du|_C$, so that
$$
\mu_C=du_1\;du_2\;du_3
$$
is the standard Lebesgue measure on the Euclidean space $C$.

Given a nonzero measure $\mu$ with compact support on a finite
dimensional real affine space $V$, the {\it barycenter} (or {\it
  centroid}) of $\mu$ is the point $\operatorname{bar}(\mu)$ of $V$
defined by
$$
\operatorname{bar}(\mu)=\frac{1}{\mu(V)}\int_{x\in V} x\;d\mu(x)\;.
$$
For instance, when $\mu$ is supported on a finite set $S$, then
$\operatorname{bar}(\mu)$ is the usual affine barycenter of the
weighted family of points $\big\{\big(s, \frac{\mu(\{s\})}
{\mu(S)}\big)\big\}_{s\in S}$.

We denote  the open ball of center $0$
  and radius $r$ in the Euclidean space $\Im\;\HH$ by  $B(r)$.
  Recall that  the radius of a finite chain $C$ is denoted by $R_C$.

\bprop Let $C$ be a chain in $\partial_\infty\hnh$ and $c\in C$.
\begin{enumerate}
\item  If $C$ is a finite chain, then the center of the
  chain $C$ is equal to the barycenter of the measure $\mu_C$:
  $$
  \cen(C)=\operatorname{bar}(\mu_C)\;.
  $$
\item If $C$ is a vertical chain, there is a unique diffeomorphism
  $\tau=\tau_{C}:\,\Im\;\HH\ra C-\{\infty\}$ such that
  $\tau^*\omega=dv$, up to precomposition of $\tau$ by a vertical
  Heisenberg translation. For every Heisenberg translation or rotation
  $\ga$, we have $\tau_{\ga C} =\ga\circ\tau_{C}$.
\item
  If $C$ is a finite chain, there exists a smooth diffeomorphism
  $\tau=\tau_{C,c}$ from $B(2\pi R_C^2)$ to $C-\{c\}$, admitting a
  continuous extension to $\partial B(2\pi R_C^2)$ sending this sphere
  to $c$, such that $\tau^*\omega=dv$.  This mapping is unique up to
  postcomposition by a Heisenberg rotation preserving $C$ and $c$, and
  $2\pi R_C^2$ is the unique radius for which such a mapping exists.
  
 ~~~ For every Heisenberg translation or rotation $\ga$, we have
 $\tau_{\ga C,\ga c} =\ga\circ\tau_{C,c}$.
\end{enumerate}
\eprop

\dem
(1) Note that $\HHeis_{4n-1}=\HH^{n-1}\times\,\Im\;\HH$ has a natural
structure of a real affine space, and that the elements of
$\operatorname{PB}_q$ act by affine transformations on
$\HHeis_{4n-1}$.  This can be seen by saying that $\HHeis_{4n-1}$,
identified with the boundary of the projective model of $\hnh$ minus
$\{\infty \}$, is a $\operatorname{PB}_q$-invariant affine subspace of
the affine chart of the quaternionic projective space defined by the
quaternionic projective hyperplane $\{[z_0:z:z_n]\in\PP^n_{\rm r}
(\HH): z_n=0\}$, and that the quaternionic projective transformations
preserving this hyperplane (as are the elements of
$\operatorname{PB}_q$) act by affine transformations on the associated
affine chart.  Another way is to check, by an easy computation, that
the Heisenberg translations, rotations and dilations preserve the
barycenters in the real affine space $\HH^{n-1}\times\,\Im\;\HH\;$:
For instance, for all $(\zeta_0,u_0),
(\zeta,u),(\zeta',u')\in\HHeis_{4n-1}$ and $t\in [0,1]$, we have
$$
(\zeta_0,u_0)\cdot\big(t \,(\zeta,u)+(1-t)(\zeta',u')\big)=
t\,(\zeta_0,u_0)\cdot(\zeta,u)+(1-t)(\zeta_0,u_0)\cdot(\zeta',u')\;.
$$
In particular, the barycenters of measures $\mu$ with compact support
on $\HHeis_{4n-1}$ are invariant under the Heisenberg translations,
rotations and dilations : For every such transformation $\ga$, we have
\begin{equation}\label{eq:equivbary}
  \operatorname{bar}(\ga_*\mu)=\ga\operatorname{bar}(\mu)\;.
\end{equation}

In order to prove Assertion (1), by Equations \eqref{eq:equivcenter}
and \eqref{eq:equivbary}, and by the transitivity properties of the
Heisenberg translations and dilations on chains, we may assume that
$n=2$ and that $C$ is a Euclidean sphere with center $(0,0)$ and
radius $1$ in the horizontal subspace $\{(\zeta,u)\in \HH^{n-1}\times
\HH :u=0\}$.  Since the $\Im\;\HH$-valued $1$-form $\omega|_C$ is
invariant under the Heisenberg rotations, the volume form $\omega_1
\wedge \omega_2 \wedge \omega_3$ on $C$ is invariant under the
Heisenberg rotations.  Since the only measure on $C$ invariant under
the Heisenberg rotations is, up to a scalar multiple, the Lebesgue
measure on the Euclidean sphere $C$, the measure $\mu_C$ is a multiple
of the Lebesgue measure on $C$.  (This can also be proved by a direct
computation: On $C$, with $\zeta=\zeta_0+\zeta_1\,i+\zeta_2\,j
+\zeta_3\,k$, we have $\omega_1\wedge\omega_2\wedge\omega_3 =-8
\sum_{i=0}^4(-1)^i\zeta_i\, d\zeta_0\wedge\cdots\widehat{d\zeta_i}
\cdots \wedge d\zeta_4$ .) Since the barycenter of this measure is
exactly the origin $(0,0)$, which is the center of the finite chain
$C$, this proves Assertion (1).

\medskip
(2) First assume that $C$ is the standard vertical chain
$$
C_\infty = \{(\zeta,u)\in\HH^{n-1}\times \,\Im\;\HH
: \zeta=0\} \cup \{\infty\}\;.
$$
Let $\tau=\tau_{C_\infty}:v\mapsto (0,v)$. Then $\tau$ is a
diffeomorphism from $\Im\;\HH$ onto $C-\{\infty\}$, such that
$\tau^*(du-2\, \Im(\, \overline{\zeta}\, d\zeta)) = dv$. For every
vertical Heisenberg translation $\ga$, the map $\ga\circ\tau$ is also
a diffeomorphism from $\Im\;\HH$ onto $C-\{\infty\}$, and since
$\omega$ is invariant under the Heisenberg translations, we also have
$(\ga\circ\tau)^*\omega = dv$.

If $\sigma:\Im\;\HH\ra C-\{\infty\}$ is another diffeomorphism such
that $\sigma^*\omega=dv$, then for every $v\in\,\Im\;\HH$, we have
$\sigma'(v)-\tau'(v)\in TC\cap \ker \omega=\{0\}$, thus the maps
$\sigma$ and $\tau$ differ by an element of the vector subspace
$C$. Therefore there exists a vertical Heisenberg translation $\ga$
such that $\sigma=\ga\circ\tau$.

Now, if $C$ is another vertical chain, there exists a composition
$\ga$ of Heisenberg translations and rotations such that $C=\ga
C_\infty$.  Defining $\tau_{C}= \ga\circ \tau_{C_\infty}$ gives a
diffeomorphism from $\Im\;\HH$ onto $C-\{\infty\}$ such that
${\tau_{C}}^*\omega=dv$, by the invariance of $\omega$ under the
Heisenberg translations and rotations. This proves Assertion (2).

\medskip
(3) First assume that $C$ is the Euclidean $3$-sphere
$$
\{(\zeta,u)\in\HH^{n-1}\times \,\Im\;\HH\; :\; \n(\zeta_1)=R^2
\;\;{\rm and}\;\;u=\zeta_2=\dots=\zeta_{n-1}=0\}\;,
$$
and that $c=(\zeta_c=(-R,0,\dots,0), u_c=0)$. Note that $R$ is the
radius of the finite chain $C$. By the properties of the exponential
map of the Lie group of unit quaternions, whose tangent space at the
identity element $1$ is $\Im\;\HH$, the smooth map
$$
\tau=\tau_{C,c}:v\mapsto (\zeta=(R\,e^{-v/(2R^2)},0,\dots,0),u=0)
$$
from $\Im\;\HH$ to $C$ is a diffeomorphism from $B(2\pi R^2)$ onto
$C-\{c\}$. It extends continuously (and even smoothly) to the sphere
$\partial B(2\pi R^2)$, mapping this sphere to $c$. Considering
$\zeta$ as a function of $v$, we have $d\zeta=(-\frac{1}{2R}\;
e^{-v/(2R^2)}dv, 0,\dots,0)$. Hence, since $v$ and $dv$ are purely
imaginary quaternions, we have
$$
\tau^*\omega=-2\;\Im(\, \overline{\zeta}\,d\zeta)=-2\;\Im\big(
\big(R\,e^{-\,\overline{v}\,/(2R^2)}\big)\big(-\frac{1}{2R}\;
e^{-v/(2R^2)}dv\big)\big)=dv\;.
$$
The uniqueness of $\tau$ up to postcomposition by a Heisenberg rotation
preserving $C$ and $c$, and the extension to the other chains, follow
as previously from the fact that the chains are transverse to the
quaternionic contact structure on $\HHeis_{4n-1}$ and by invariance of
the calibration $\omega$ under the Heisenberg translations and
rotations.  \cqfd

\section{Counting and equidistribution of arithmetic chains \\
  in hyperspherical geometry}
\label{sect:countequidistrib}

In this section, we prove (generalised versions of) Theorems
\ref{theo:countingintro} and \ref{theo:equidistribintro} of the
introduction. We start by recalling a general statement, coming from a
special case of the main results of \cite{ParPau17ETDS}, that has been
explicited in \cite{ParPau20a}.

Let $\Ga$ be a lattice in $\PU_q$. Let $D^-$ and $D^+$ be nonempty
proper closed convex subsets of $\hnh$, with stabilisers $\Ga_{D^-}$
and $\Ga_{D^+}$ in $\Ga$ respectively, such that the families $(\ga
D^-)_{\ga\in\Ga/\Ga_{D^-}}$ and $(\ga D^+)_{\ga\in\Ga/\Ga_{D^+}}$ are
locally finite in $\hnh$.
For all $\ga,\ga'$ in $\Ga$, the convex sets $\ga D^-$ and $\ga' D^+$
have a common perpendicular if and only if their closures
$\overline{\ga D^-}$ and $\overline{\ga' D^+}$ in $\hnh\cup
\partial_\infty \hnh$ do not intersect.  We denote by
$\alpha_{\ga,\,\ga'}$ this common perpendicular, starting from $\ga
D^-$ at time $t=0$, and by $\ell(\alpha_{\ga,\,\ga'})$ its length.
The {\em multiplicity} of $\alpha_{\ga,\ga'}$ is
$$
m_{\ga,\ga'}=
\frac 1{\card(\ga\Ga_{D^-}\ga^{-1}\cap\ga'\Ga_{D^+}{\ga'}^{-1})}\,,
$$ 
which equals $1$ for all $\ga,\ga'\in\Ga$ when $\Ga$ acts freely on
$T^1\hnh$ (for instance when $\Ga$ is torsion-free). For all $s> 0$
and $x\in\partial D^-$, let
$$
m_s(x)=\sum_{\ga\in \Ga/\Ga_{D^+}\;:\;
  \overline{D^-}\,\cap \,\overline{\ga D^+}\,= \emptyset,\;
\alpha_{e,\, \ga}(0)=x,\; \ell(\alpha_{e,\, \ga})\leq s} m_{e,\ga}
$$ 
be the multiplicity of $x$ as the origin of common perpendiculars
with length at most $t$ from $D^-$ to the elements of the $\Ga$-orbit
of $D^+$. 

For every $s> 0$, let
$$
\N_{D^-,\,D^+}(s)=
\sum_{
(\ga,\,\ga')\in \Ga\bs((\Ga/\Ga_{D^-})\times (\Ga/\Ga_{D^+}))\;:\;
\overline{\ga D^-}\,\cap \,\overline{\ga' D^+}\,=\emptyset,\; 
\ell(\alpha_{\ga,\, \ga'})\leq s} m_{\ga,\ga'}
\;,
$$ 
where $\Ga$ acts diagonally on $\Ga\times\Ga$. When $\Ga$ has no
torsion, $\N_{D^-,\,D^+}(s)$ is the number (with multiplicities coming
from the fact that $\Ga_{D^\pm}\bs D^\pm$ is not assumed to be
embedded in $\Ga\bs\hnh$) of the common perpendiculars of length at
most $s$ between the images of $D^-$ and $D^+$ in $\Ga\bs\hnh$.

The following statement is a special case of
\cite[Thm.~8.1]{ParPau20a}.  We denote by $\Delta_x$ the unit Dirac
mass at a point $x$.

\btheo\label{theo:quaterhyperbo} Let $D^-$ be a horoball in $\hnh$
centred at a parabolic fixed point of $\Ga$ and let $D^+$ be a
quaternionic geodesic line in $\hnh$ such that $\Ga_{D^+}\bs D^+$ has
finite volume.  Let $m^+$ be the order of the pointwise stabiliser of
$D^+$ in $\Ga$ and let
$$
c(D^-,D^+)=\frac{2\,(n-1)\,(2n-1)}{\pi^2\,m^+}
\frac{\Vol(\Ga_{D^-}\bs D^-)\Vol(\Ga_{D^+}\bs D^+)}
{\Vol(\Ga\bs\hnh)}\,. 
$$
There exists $\kappa>0$ such that, as
$s\ra+\infty$,
$$
\N_{D^-,\,D^+}(s)=c(D^-,D^+)\;
e^{(4n+2)\,s}\;\big(1+\operatorname{O}(e^{-\kappa s})\big)\;.
$$
Furthermore, the origins of the common perpendiculars from $D^-$ to
the images of $D^+$ under the elements of $\Ga$ equidistribute in
$\partial D^-$ to the induced Riemannian measure: as $s\ra+\infty$,
\begin{equation}\label{eq:distribhorobcomplexhyp}
\frac{2\,(2n+1)\,\Vol(\Ga_{D^-}\bs D^-)}{c(D^-,D^+)}
\;e^{-(4n+2)\,s}\;\sum_{x\in\partial D^-} m_{s}(x) \;
\Delta_{x}\;\weakstar\; \vol_{\partial D^-}\,.
\;\;\;\Box
\end{equation}
\etheo

For smooth functions $\psi$ with compact support on $\partial D^-$,
there is an error term in the equidistribution claim of Theorem
\ref{theo:quaterhyperbo} when the measures on both sides are evaluated
on $\psi$, of the form $\bigO(e^{-\kappa s}\,\|\psi\|_\ell)$ where
$\kappa>0$ and $\|\psi\|_\ell$ is the Sobolev norm of $\psi$ for some
$\ell\in\NN$.

\medskip
From now on, we assume that $n=2$. Let $A$, $D_A$, $m_A$ and $\OOO$ be
as in the Introduction. We denote by $|\,\OOO^\times|$ the order of
the unit group of $\OOO$, equal to $24$ if $D_A=2$, to $12$ if
$D_A=3$, or else to $2$, $4$ or $6$. See for instance
\cite{Vigneras80}.  As usual, by $\prod_{p|D_A}$, we mean a product
where $p$ ranges over the prime positive numbers dividing $D_A$.

For every chain $C$ in $\partial_\infty\hdh$, let $L_C$ be the
quaternionic projective line in $\PP^2_{\rm r}(\HH)$ such that
$C=L_C\cap \partial_\infty\hdh$, and let $D_C=L\cap\hdh$ be the
associated quaternionic geodesic line. For every finite index subgroup
$G$ of the arithmetic lattice $\PU_q(\OOO)$, we denote by $G_C$ the
stabiliser of $C$ in $G$, by $G_\infty$ the stabiliser of $\infty$ in
$G$, and by $\covol_G(C)$ the volume of the orbifold $G_C\bs D_C$ for
the Riemannian metric of constant sectional curvature $-1$ on the real
hyperbolic $4$-space $D_C$.  Recall that a chain $C$ is arithmetic
over $\OOO$ if and only if the stabiliser in $\PU_q(\OOO)$ (or
equivalently in $G$) of the quaternionic geodesic line $D_C$ has
finite covolume on $D_C$.

\btheo \label{theo:countchain} Let $C_0$ be an arithmetic chain over a
maximal order $\OOO$ in a definite quaternion algebra over $\QQ$. Let
$G$ be a finite index subgroup of $\PU_q(\OOO)$. Then there exists a
constant $\kappa>0$ such that, as $\epsilon>0$ tends to $0$, the
number $\psi_{C_0,\,G}(\epsilon)$ of chains modulo $G_\infty$ in the
$G$-orbit of $C_0$ with $d_{\rm Cyg}$-diameter at least $\epsilon$ is
equal to
$$
\frac{25515\;2^{23}\;D_A^2\;\covol_G(C_0)\,[\PU_q(\OOO)_\infty:G_\infty]}
{\pi^6\;m_{C_0,G}\;m_A\;|\OOO^\times|^2\;\prod_{p|D_A}(p-1)(p^2+1)(p^3-1)\;
[\PU_q(\OOO):G]}
\;\epsilon^{-10}\big(1+\bigO(\epsilon^\kappa)\big)\;,
$$
where $m_{C_0,G}$ is the order of the pointwise stabiliser of
$D_{C_0}$ in $G$.
\etheo

Recall that the center $\operatorname{cen}(C)$ of a finite chain $C$
is the image of $\infty=[1:0:0]$ under the reflexion on $L_C$. The
following result is an equidistribution result in the quaternionic
Heisenberg group of the centers of the arithmetic chains in a given
orbit under (a finite index subgroup of) $\PU_q(\OOO)$.

\btheo\label{theo:equidischain} Let $C_0$ and $G$ be as in Theorem
\ref{theo:countchain}. As $\epsilon>0$ tends to $0$, we have
\begin{multline*}
\frac{m_{C_0,G}\,m_A\,\pi^6\,\;\prod_{p|D_A}(p-1)(p^2+1)(p^3-1)\;
[\PU_q(\OOO):G]}
{25515\;2^{24}\covol_G(C_0)} \;\epsilon^{10}\\
\sum_{\scriptsize\begin{array}{c}C\in  G\cdot C_0\\
\diam_{d_{\rm Cyg}}(C)\geq \epsilon\end{array}}\; 
\Delta_{\operatorname{cen}(C)}\;\weakstar\;\haarhheis\,.
\end{multline*}
\etheo

As in Theorem \ref{theo:quaterhyperbo}, there exist $\kappa>0$ and
$\ell\in\NN$ such that for every smooth function $\psi$ with compact
support on $\HHeis_7$, there is an error term in this equidistribution
result when the measures on both sides are evaluated on $\psi$, of the
form $\bigO(s^{-\kappa}\,\|\psi\|_\ell)$ where $\|\psi\|_\ell$ is the
Sobolev norm of $\psi$.

\medskip
We begin by a technical result used in the proofs of the above
theorems, which does not require the assumption $n=2$. Recall that
$d''_{\rm Cyg}$ is the modified Cygan distance defined in Section
\ref{sect:quathypheis}.

\blemm\label{lem:diametre}
For every $m$-chain $C$ in $\hnh$, we have $\diam_{d_{\rm Cyg}}(C)=
\sqrt{2}\;\diam_{d''_{\rm Cyg}}(C)$.
\elemm

\dem
If $C$ is a vertical $m$-chain, then both diameters are $+\infty$.  We
hence assume that $C$ is finite. Since the Heisenberg translations and
rotations preserve $d_{\rm Cyg}$ and $d''_{\rm Cyg}$, and by the
transitivity properties of the Heisenberg translations and rotations
on the set of $m$-chains (see Section \ref{subsec:transitivity}), we
may assume that $C$ is a Euclidean sphere centered at $(0,0)$ with
dimension $4m-1$, contained in the horizontal plane $\HH^{n-1}\times
\{0\}$ of $\HHeis_{4n-1}$.

On this horizontal plane, the modified Cygan distance coincides with
$\frac{1}{\sqrt{2}}$ times the Cygan distance, by Equations
\eqref{eq:defidistCyg} and \eqref{eq:defidistCyg''}. Hence the
diameter of the $m$-chain $C$ for the modified Cygan distance is equal
to $\frac{1}{\sqrt{2}}$ times its diameter for the Cygan distance.
This proves the result.
\cqfd

\bigskip
\noindent{\bf Proof of Theorem \ref{theo:countchain} and Theorem
  \ref{theo:equidischain}. }  The diameter of a chain for the Cygan
distance is invariant under the stabiliser in $\PU_q$ of the
horosphere $\partial \H_1$, hence is invariant under $G_\infty$. The
counting function $\psi_{C_0,\,G}$ is thus well defined.

Note that $\H_1$ is a horoball centered at the fixed point of a
parabolic element in $\PU_q(\OOO)$ (take the vertical Heisenberg
translation by $(0,2u)$ for any nonzero $u\in \OOO\cap\Im\,\HH$~). We will
apply Theorem \ref{theo:quaterhyperbo} with $\Ga=G$, with $D^-=\H_1$,
which is hence a horoball centered at the fixed point of a parabolic
element in $G$, and with $D^+= D_{C_0}$, which is the quaternionic
geodesic line in $\hdh$ with boundary at infinity equal to $C_0$. In
particular $m^+=m_{C_0,G}$.

Let us compute the constant $c(D^-,D^+)$ appearing in the statement of
Theorem \ref{theo:quaterhyperbo}. We have $\Vol(\Ga\bs\hdh)=
[\PU_q(\OOO):G]\Vol(\PU_q(\OOO)\bs \hdh)$, where, by
\cite[Thm.~1.4]{ParPau20a},
$$
\Vol(\PU_q(\OOO)\bs \HH^2_\HH)=
\frac{\pi^4\;m_A}{42525\;2^{13}}\;\prod_{p|D_A}(p-1)(p^2+1)(p^3-1)\;,
$$
and  by \cite[Lem.~8.4]{ParPau20a},
\begin{equation}\label{eq:volcusp}
\Vol(\Ga_{D^-}\bs D^-)=[\PU_q(\OOO)_\infty:G_\infty]
\Vol(\PU_q(\OOO)_{\H_1}\bs\,\H_1)=
\frac{D_A^2\,[\PU_q(\OOO)_\infty:G_\infty]}{160\,|\,\OOO^\times|^2}\;.
\end{equation}
By definition, we have
$$
\Vol(\Ga_{D^+}\bs D^+)=16\;\covol_G(C_0)\,,
$$
since the sectional curvature of $D^+$ is constant $-4$ and $D^+$
has real dimension $4$. We
hence have
\begin{equation}\label{eq:cee}
c(D^-,D^+)=
\frac{25515\;2^{13}\;D_A^2\;\covol_G(C_0)\,[\PU_q(\OOO)_\infty:G_\infty]}
{\pi^6\;m_{C_0,G}\;m_A\;|\OOO^\times|^2\;
\prod_{p|D_A}(p-1)(p^2+1)(p^3-1)\;[\PU_q(\OOO):G]}\,. 
\end{equation} 

\medskip 
Let $g\in G$ be such that the quaternionic geodesic line $gD^+$ is
disjoint from $\H_1$ (which is the case except for $g$ in finitely
many double classes in $G_{\H_1}\bs G/G_{D^+}$). Let $\delta_g$ be the
common perpendicular from $\H_1$ to $gD^+$. Its length
$\ell(\delta_g)$ is the minimum of the distances from $\H_1$ to a
geodesic line between two points of $\partial_\infty(gD^+)=g
C_0$. Hence, by Lemmas \ref{lem:distgeod} and \ref{lem:diametre}, we
have
\begin{align}\label{eq:relatlongcomperpchain}
\ell(\delta_g)= &
\min_{x,y\in gC_0,\;x\neq y} d(\H_1,\;]x,y[\,) 
=-\max_{x,y\in gC_0,\;x\neq y}\ln \frac{d''_{\rm Cyg}(x,y)}{\sqrt 2}
\nonumber\\
= & 
-\ln \frac{\diam_{d''_{\rm Cyg}}(gC_0)}{\sqrt{2}}=
-\ln \frac{\diam_{d_{\rm Cyg}}(gC_0)}{2}\,.
\end{align}

Respectively by the definition of the counting function
$\psi_{C_0,\,G}$ in the statement of Theorem \ref{theo:countchain},
since the stabiliser of $C_0$ in $G$ is equal to $G_{D^+}$, by
Equation \eqref{eq:relatlongcomperpchain}, by Theorem
\ref{theo:quaterhyperbo}, and by Equation \eqref{eq:cee}, we have, as
$\epsilon>0$ tends to $0$,
\begin{align*}
&\psi_{C_0,\,G}(\epsilon)\\
=\;&\card\;\; G_\infty\bs\{C\in  G\cdot C_0\;:\; 
\diam_{d_{\rm Cyg}}(C)\geq \epsilon\}\\  =\;&
\card\{[g]\in G_\infty\bs G/G_{D_{C_0}}\;:\; 
\diam_{d_{\rm Cyg}}(gC_0)\geq \epsilon\}\\ =\; & 
\card\{[g]\in G_{\H_1}\bs G/G_{D_{C_0}}\;:\; \ell(\delta_g)
\leq -\ln\frac{\epsilon}{2}\} +\bigO(1)\\ =\;&
\N_{D^-,\,D^+}(-\ln\frac{\epsilon}{2})+\bigO(1)=
c(D^-,D^+)\;e^{-10\ln\frac{\epsilon}{2}}
\big(1+\bigO(e^{\kappa \ln\frac{\epsilon}{2}})\big)\\ =\;&
\frac{25515\;2^{23}\;D_A^2\;\covol_G(C_0)\,[\PU_q(\OOO)_\infty:G_\infty]}
{\pi^6\;m_{C_0,G}\;m_A\;|\OOO^\times|^2\;\prod_{p|D_A}(p-1)(p^2+1)(p^3-1)\;
[\PU_q(\OOO):G]}
\;\epsilon^{-10}\big(1+\bigO(\epsilon^\kappa)\big)\;.
\end{align*}

This proves Theorem \ref{theo:countchain}. Let us now prove Theorem
\ref{theo:equidischain}.

\medskip We apply the equidistribution result in Equation
\eqref{eq:distribhorobcomplexhyp} of the origins
$\operatorname{or}(\delta_g)$ of the common perpendiculars $\delta_g$
from $D^-=\H_1$ to the images $gD^+$ for $g\in G$. As $s\ra+\infty$,
we hence have, using Equations \eqref{eq:cee} and \eqref{eq:volcusp},
\begin{multline}
\frac{m_{C_0,G}\,m_A\,\pi^6\,\;\prod_{p|D_A}(p-1)(p^2+1)(p^3-1)\;
[\PU_q(\OOO):G]}{25515\;2^{17}\covol_G(C_0)}
\;e^{-10\,s}\\ \sum_{[g]\in G/G_{D^+}\,:\; \ell(\delta_g)\leq s}\;
\Delta_{\operatorname{or}(\delta_g)}\;\weakstar\; \vol_{\partial \H_1}\,.
\label{eq:distribcentchain}
\end{multline}

Let $f:\partial_\infty\hdh-\{\infty\}=\HHeis_7\ra \partial \H_1$ be
the orthogonal projection map, which is the homeomorphism $(w_0,w)
\mapsto (w_0+\frac{1}{2},w)$. The pushforward of the Haar measure
$\haarhheis$ by $f$ is
\begin{equation}\label{eq:volheis}
f_*\haarhheis=8 \vol_{\partial \H_1}\,,
\end{equation}
see for example the end of the proof of Theorem 8.3 in
\cite{ParPau20a}.

Note that, for every chain $C$, if $r_C$ is the reflexion on the
quaternionic projective line containing $C$, then the geodesic line
from $\infty$ to $\operatorname{cen}(C)=r_C(\infty)$, being invariant
under $r_C$, is orthogonal to the quaternionic geodesic line with
boundary at infinity $C$.  Hence for every $g\in G$, we have
$$
f^{-1}(\operatorname{or}(\delta_g))=\operatorname{cen}(gC_0)\;.
$$

Let us use in Equation \eqref{eq:distribcentchain} the change of
variables $s=-\ln\,\frac{\epsilon}{2}$ and the continuity of the
pushforward of measures by $f^{-1}$. By Equations
\eqref{eq:relatlongcomperpchain} and \eqref{eq:volheis}, as
$\epsilon>0$ tends to $0$, we obtain that the measures
$$
\frac{m_{C_0,G}\,m_A\,\pi^6\,\;\prod_{p|D_A}(p-1)(p^2+1)(p^3-1)\;
[\PU_q(\OOO):G]}{25515\;2^{24}\covol_G(C_0)}
\;\epsilon^{10}\!\!\!
\sum_{\scriptsize\begin{array}{c}[g]\in  G/G_{D^+}\\
\diam_{d_{\rm Cyg}}(gC_0)\geq \epsilon\end{array}}
\Delta_{\operatorname{cen}(gC_0)}
$$
weak-star converge to the Haar measure $\haarhheis$.  This proves
Theorem \ref{theo:equidischain}.
\cqfd

\bigskip
\noindent{\bf Example. } Let $C_0=\big\{[w_0:0:1]\in\PP^2_{\rm r}
(\HH):\tr w_0=0\big\}$ be the standard vertical chain in
$\partial_\infty\hdh$, which is the intersection of
$\partial_\infty\hdh$ with the quaternionic projective line
$L_{C_0}=\big\{[z_0:z_1:z_2]\in\PP^2_{\rm r} (\HH): z_1=0\big\}$.

An element $\pm\begin{pmatrix} a & \gamma^* & b\\ \alpha & A & \beta
\\ c & \delta^* & d\end{pmatrix}$ of $\PU_q$ preserving the
quaternionic geodesic line $L_{C_0}\cap\hdh$ satisfies $\alpha
w_0+\beta=0$ for all $w_0\in\HH$ with $\tr w_0>0$.  Thus, $\alpha=
\beta=0$, and Equations \eqref{eq:equationsUq} (or rather the similar
equations obtained by the formula $X X^*=I_{n+1}$ instead of $X^* X=
I_{n+1}$) imply that $\gamma= \delta=0$.  Using again Equations
\eqref{eq:equationsUq}, we see that the stabiliser of $L_{C_0}$
consists of the elements $\begin{pmatrix} a & 0 & b\\ 0 & A & 0 \\ c &
0 & d\end{pmatrix}$ such that $\tr(\,\overline{c}\,a)=
\tr(\,\overline{d}\,b) =0$, ~$\overline{c}\,b+\overline{a}\,d=1$ and
$A\in\OOO^\times$.  Thus,
$$
\covol_{\PU_q(\OOO)}(C_0)=\frac{\pi^2}{1080}\prod_{p|D_A}(p-1)(p^2+1)
$$
by \cite[Thm.~2.5]{BreHel96}.

The pointwise stabiliser of $C_0$ in $\PU_q(\OOO)$ consists of the
diagonal elements with $a=d=\pm1$ and $A\in \OOO^\times$, giving
$m_{C_0,\PU_q(\OOO)}=|\OOO^\times|$.

Theorems \ref{theo:countchain} and \ref{theo:equidischain} then give
$$
\psi_{C_0,\,\PU_q(\OOO)}(\epsilon)=
\frac{189\;2^{20}\;D_A^2}{\pi^4\;m_A\;|\OOO^\times|^3\;\prod_{p|D_A}(p^3-1)}
\;\epsilon^{-10}\big(1+\bigO(\epsilon^\kappa)\big)\;, 
$$
and
$$
\frac{\pi^4\,m_A\,|\OOO^\times|\;\prod_{p|D_A}(p^3-1)}{189\;2^{21}}
\;\epsilon^{10}
\sum_{C\in \PU_q(\OOO)\cdot C_0\;:\;\diam_{d_{\rm Cyg}}\;C\geq \epsilon}\;
\Delta_{\operatorname{cen}(C)}\;\weakstar\;\haarhheis\,.
$$

{\small \bibliography{../biblio} }

\bigskip
{\small
\noindent \begin{tabular}{l} 
Department of Mathematics and Statistics, P.O. Box 35\\ 
40014 University of Jyv\"askyl\"a, FINLAND.\\
{\it e-mail: jouni.t.parkkonen@jyu.fi}
\end{tabular}
\medskip

\noindent \begin{tabular}{l}
Laboratoire de math\'ematique d'Orsay, UMR 8628 CNRS\\
Universit\'e Paris-Saclay,
91405 ORSAY Cedex, FRANCE\\
{\it e-mail: frederic.paulin@universite-paris-saclay.fr}
\end{tabular}
}

\end{document}